# $O(N)$ Iterative and $O(NlogN)$ Fast Direct Volume Integral Equation Solvers with a Minimal-Rank $\mathcal{H}^2$-Representation for Large-Scale 3-D Electrodynamic Analysis

Saad Omar, *Member, IEEE* and Dan Jiao, *Fellow, IEEE*

*Abstract*—Linear complexity iterative and log-linear complexity direct solvers are developed for the volume integral equation (VIE) based general large-scale electrodynamic analysis. The dense VIE system matrix is first represented by a new cluster-based multilevel low-rank representation. In this representation, all the admissible blocks associated with a single cluster are grouped together and represented by a single low-rank block, whose rank is minimized based on prescribed accuracy. From such an initial representation, an efficient algorithm is developed to generate a minimal-rank $\mathcal{H}^2$-matrix representation. This representation facilitates faster computation, and ensures the same minimal rank's growth rate with electrical size as evaluated from singular value decomposition. Taking into account the rank's growth with electrical size, we develop linear-complexity $\mathcal{H}^2$-matrix-based storage and matrix-vector multiplication, and thereby an $O(N)$ iterative VIE solver regardless of electrical size. Moreover, we develop an $O(NlogN)$ matrix inversion, and hence a fast $O(NlogN)$ *direct* VIE solver for large-scale electrodynamic analysis. Both theoretical analysis and numerical simulations of large-scale 1-, 2- and 3-D structures on a single-core CPU, resulting in millions of unknowns, have demonstrated the low complexity and superior performance of the proposed VIE electrodynamic solvers.

*Index Terms*—Volume integral equations, fast direct solvers, fast solvers, electrodynamic, linear complexity solvers, scattering, radiation, electromagnetic modeling, three dimensional structures

## I. INTRODUCTION

The volume integral equation (VIE) based methods [1]–[6] offer great flexibility in modeling both complicated geometries and inhomogeneous materials in open-region settings. From a computational perspective, integral equation methods lead to dense matrices. The size of these dense matrices, for volume based analysis, increases cubically with the size of the objects under study. Therefore, the advantages of the VIE-based analyses can be fully accentuated only if they can be performed with low computational complexity.

Existing fast VIE solvers for solving large-scale electrodynamic problems are, in general, iterative solvers. These include

Manuscript received on June 22, 2016. This work was supported by a grant from NSF under award No. 1619062, a grant from SRC (Task 1292.073), and a grant from DARPA under award HR0011-14-1-0057. The authors are with the School of Electrical and Computer Engineering, Purdue University, West Lafayette, IN 47907 USA.

methods like Fast Fourier Transform (FFT) [7]–[13], low-rank compression [14], [15], Fast Multipole Method (FMM) [16]–[20], and others. The memory requirements in such methods scale at best as $O(N)$, whereas a single matrix-vector multiplication can cost as small as $O(NlogN)$ in time complexity, where $N$ is matrix size. The overall complexity for such iterative solvers is $O(N_{rhs}N_{it}NlogN)$, where $N_{rhs}$ is the number of right hand sides, and $N_{it}$ is the number of iterations. When the number of right hand sides under analysis and/or the number of iterations are large, iterative solvers become inefficient since an entire iteration procedure has to be repeated for each right hand side.

There have been significant contributions in fast direct solvers [23]–[41], for the analysis of problems ranging from circuits to scattering problems. For VIE analysis, although $O(N)$ direct VIE solvers have been developed for full-wave general 3-D circuit analysis [40], no $O(Nlog^\alpha N)(\alpha \geq 0)$ fast solvers have yet been made possible for electrically large analysis. The main contribution of this paper is such a fast direct VIE solver whose inversion complexity is $O(NlogN)$ irrespective of electrical size, in addition to a fast iterative VIE solver whose matrix-vector multiplication complexity is $O(N)$.

To achieve these low complexities, one of the key challenges is to compactly represent the dense VIE system matrix into a reduced set of parameters, despite the large and electrical-size dependent rank in the off-diagonal blocks. For example, if one uses a degenerate approximation of IE kernels, i.e., separating sources from observers in approximating Green's function like that in the FMM or an $\mathcal{H}^2$-based circuit solver [40], he would obtain an asymptotically full-rank representation of the original electrodynamic IE operator. This challenge is overcome in this work by finding a minimal-rank representation to approximate the VIE operator based on prescribed accuracy, which does not separate sources from observers. The minimal rank required for representing the VIE operator will be analyzed in Section II. Such a minimal-rank representation can be obtained from singular value decomposition (SVD). However, a brute-force SVD is computationally expensive. We hence develop an efficient algorithm to represent the original dense system matrix using its minimal rank required by accuracy, while avoiding the huge computational cost of SVD. This algorithm will be detailed in the following Section III. In Section IV, we present proposed fast $O(N)$ iterative and $O(NlogN)$ direct inverse VIE solvers for electrodynamic analysis, enabled by the proposed minimal-rank as well as nested representation. In Section V, numerical results are presented to demonstrate the performance of the proposed iterative and direct solvers for arbitrary electrodynamic analysis. Section VI relates to our conclusions. This paper is a significant expansion of our conference paper [54] from theory, algorithm, and numerical experiments perspectives. The algorithms developed in this work for obtaining a minimal-rank $\mathcal{H}^2$-representation of complex-valued dense matrices are purely algebraic and kernel independent. In addition to VIE, they can also be applied to other IE operators.

It is worth mentioning that the algorithm developed in this work for obtaining a minimal-rank $\mathcal{H}^2$-matrix is very different from that in [34]. In [34], because the problem being considered is a circuit extraction problem whose electric size is small, an interpolation-based $\mathcal{H}^2$-representation is first obtained, which is then converted to a minimal-rank $\mathcal{H}^2$-matrix. Hence, the algorithm developed therein is to convert an initial $\mathcal{H}^2$-matrix whose rank is not minimal to a new $\mathcal{H}^2$ matrix, whose rank is minimal for a prescribed accuracy. In this work, the interpolation-based $\mathcal{H}^2$-representation cannot be used since it upfront would produce a full-rank representation when handling electrically large Green's function. Therefore, our algorithm in this work for obtaining a minimal-rank $\mathcal{H}^2$ for an electrically large kernel is very different from the algorithm in [34], the details of which can be found from Section III. In addition, all the algorithms developed in this work are for complex-valued numerical systems, unlike the real-valued system concerned in an electrically small analysis. Our previous direct-solver work reported in [29]–[32], [40] are all based on an interpolation-based method to obtain an $\mathcal{H}^2$-matrix, which is not amenable for handling the rank's growth with electrical size in electrically large analysis. Therefore, the problem studied in this paper has not been addressed by our previous work.

## II. On the VIE Operator and Its Rank

### A. VIE Formulation for Wave-based Analysis

Consider an arbitrarily shaped 3-D inhomogeneous dielectric body of complex permittivity $\epsilon(\mathbf{r})$ occupying volume $V$, which is exposed to an incident field $\mathbf{E}^i(\mathbf{r})$.

The scattered field due to the equivalent volume polarization current $\mathbf{J}$ contributes to the total field at any point $\mathbf{r}$ in the sense as expressed in the form of the following volume integral equation,

$$\mathbf{E}^i(\mathbf{r}) = \frac{\mathbf{D}(\mathbf{r})}{\epsilon(\mathbf{r})} -$$
$$\int_V \left[ \mu_0 \omega^2 \kappa(\mathbf{r}') \mathbf{D}(\mathbf{r}') + \nabla \left( \nabla' \cdot (\kappa(\mathbf{r}) \frac{\mathbf{D}(\mathbf{r}')}{\epsilon_0}) \right) \right] g(\mathbf{r}, \mathbf{r}') dv', \quad (1)$$

where $g(\mathbf{r}, \mathbf{r}') = e^{-jk_0|\mathbf{r}-\mathbf{r}'|}/4\pi|\mathbf{r}-\mathbf{r}'|$, $\omega$ being the angular frequency, $\kappa$ the contrast ratio defined as $(\epsilon(\mathbf{r}) - \epsilon_0)/\epsilon(\mathbf{r})$, $\mathbf{D}(\mathbf{r}')$ the electric flux density, while $k_0$ is the free space wave number.

Tetrahedral based discretization is used to model the arbitrarily shaped 3-D inhomogeneous dielectric scattering body and divergence conforming Schaubert Wilson Glisson (SWG) basis functions are used [1] to expand the unknown electric flux density $\mathbf{D}(\mathbf{r}')$. Each of the SWG basis function is defined for a face of a tetrahedron.

By expanding the unknown electric flux density $\mathbf{D}(\mathbf{r}')$ in terms of SWG basis functions $\mathbf{D}_n(\mathbf{r}')$ each with a coefficient $D_n$, and then testing the resulting equation using Galerkin method with $\mathbf{D}_m(\mathbf{r})$, we obtain the following linear system of VIE,

$$\mathbf{S}D = E \quad (2)$$

where,

$$E_m = \int_{V_m} \vec{E}^i \cdot \vec{D}_m(\vec{r}) dv$$

$$\mathbf{S}_{mn} = \int_{V_m} \frac{\vec{D}_m(\vec{r})}{\bar{\epsilon}_n(\vec{r})} \cdot \vec{D}_n(\vec{r}) dv -$$
$$\mu\omega^2 \int_{V_m} \int_{V_n} \kappa_n(\vec{r}') \vec{D}_m(\vec{r}) \cdot \vec{D}_n(\vec{r}) g(\vec{r},\vec{r}') dv' dv -$$
$$\frac{1}{\epsilon_0} \Big( \int_{S_m} \int_{V_n} (\vec{D}_m(\vec{r}) \cdot \hat{n})(\nabla' \cdot \vec{D}_n(\vec{r}')) g(\vec{r},\vec{r}') dv' ds +$$
$$\int_{S_m} \int_{S_n} (\vec{D}_m(\vec{r}) \cdot \hat{n})(\nabla' \kappa_n(\vec{r}')) g(\vec{r},\vec{r}') ds' dn' ds -$$
$$\int_{V_m} \int_{V_n} (\nabla \cdot \vec{D}_m(\vec{r}))(\nabla' \cdot \vec{D}_n(\vec{r}')) g(\vec{r},\vec{r}') dv' dv -$$
$$\int_{V_m} \int_{S_n} (\nabla \cdot \vec{D}_m(\vec{r}))(\nabla' \kappa_n(\vec{r}')) g(\vec{r},\vec{r}') ds' dn' dv \Big). \quad (3)$$

As evident, each system matrix entry involves all four possible combinations of volume and surface integrals with different terms for different observation and source locations.

### B. Rank of the Electrodynamic VIE Operator

Unlike static problems, the rank of an electrodynamic IE kernel increases with electrical size for achieving a prescribed accuracy. Therefore, a fast solver built upon the low-rank property would have a higher computational complexity for solving electrically large problems as compared to electrically small problems, if no advanced algorithms are developed to effectively manage the rank's growth with electrical size. The true indicator of the rank's growth is singular value decomposition (SVD), since its resultant representation constitutes a minimal rank representation of a matrix for any prescribed accuracy. The SVD does not separate sources from observers in approximating Green's function, and it finds a minimal rank representation of the IE kernel as a whole. The SVD is computationally $O(N^3)$, and hence not practically feasible for studying the rank of electrically large IE operators. In view of the pivotal importance of this subject, a theoretical study has been carried out on the rank's growth with electrical size in integral equations [51]. A closed-form analytical expression of the rank of the coupling Green's function is derived, which has the same scaling as that depicted by SVD-based rank revealing. The findings on the rank-study are summarized as follows:

1) The rank ($k$) of the off-diagonal block, irrespective of the electrical size, is far less than the size of the block, thus the off-diagonal block has a low rank representation, i.e. $k \ll N$.

2) For static and one-dimensional configurations of sources and observers, the rank required by a prescribed accuracy remains constant irrespective of the problem size.

3) For 2- and 3-D configurations, the rank varies as square root of logarithm and linearly with the electrical size, respectively.

The findings in [51] are, in fact, consistent with the analysis in the well-known FMM-based method. As shown in [16], [52], the number of spherical harmonics required to represent Green's function for a prescribed accuracy scales *linearly* with the electric size for general 3-D problems, which agrees with the findings of [51]. This fact has actually resulted in the reduced complexity of $O(N log N)$ of an FMM-based method for one matrix-vector multiplication for solving electrically large surface IEs. Despite the reduced number of harmonics, the final representation of the dense system matrix from an FMM-based method is full rank asymptotically, which might have misled people to consider that a high-frequency kernel is not low-rank. This full-rank model, in fact, is due to a source-observer separated representation used in the FMM. When the sources are separated from observers, the Green's function, which originally only depends on the distance between sources and observers, becomes a function of complete coordinates of sources and observers. As a result, along every direction, one has to capture the linearly growing number of harmonics, leading to a full-rank representation. Such a representation is not a minimal-rank representation that does not separate sources from observers.

To numerically verify the findings related to the rank in the VIE setting concerned in this work, in Fig. 1, we plot the number of singular values and hence the rank required to maintain an accuracy of $10^{-5}$ in a 1-D type VIE configuration. The size of the two dielectric rods is kept to be 1 m geometrically separated by a distance of 2 m, while the frequency (hence, discretization also scales accordingly) sweeps to give the electric sizes as plotted on the horizontal axis. We can observe a constant rank throughout even if the electric size grows to as large as 100 wavelengths. In Fig. 2 and 3, we plot the number of singular values and hence the rank required to maintain an accuracy of $10^{-5}$ in a 2-D, and 3-D type VIE configuration, respectively. The size of the two corresponding structures is again kept to be 1 m in each dimension, geometrically separated by an $x$-directed distance of 2 m, while the frequency (hence, discretization, along each dimension, also scales accordingly) sweeps to give the electric sizes as plotted on the horizontal axes. In both configurations, we can observe that the rank's growth rate is no greater than linear with electrical size. The above numerical results obtained from the VIE operator agree very well with the theoretical findings in [51].

Since the number of unknowns in a VIE-based analysis scales with electric size in a cubic way, an error-bounded low-rank representation exists for VIE operators, irrespective of the electric size and object shape. Since SVD is computationally expensive, in the following section, we present an efficient algorithm for generating a minimal-rank $\mathcal{H}^2$-representation of the VIE dense system matrix, which has the same rank's growth rate with electrical size as that dictated by SVD.

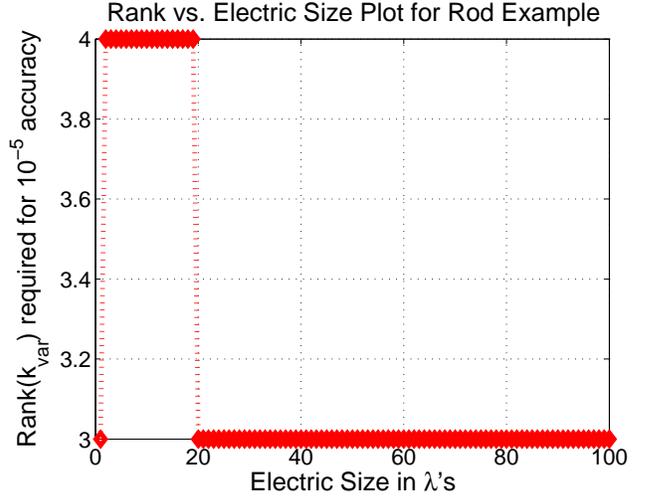

Fig. 1. Dielectric rod (1-D structure): SVD rank growth w.r.t. the electric size in an off-diagonal block.

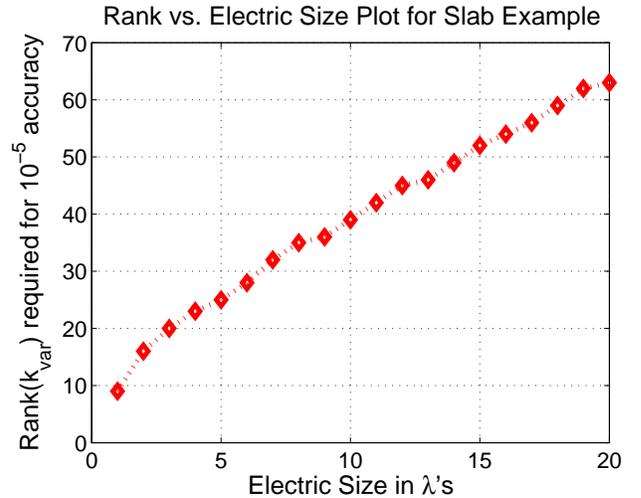

Fig. 2. Dielectric slab (2-D structure): SVD rank growth w.r.t. the electric size in an off-diagonal block.

## III. PROPOSED ALGORITHM FOR GENERATING MINIMAL-RANK $\mathcal{H}^2$-REPRESENTATION OF THE ELECTRODYNAMIC VIE OPERATOR

Starting from the original VIE operator, we develop a two-stage algorithm to obtain a minimal-rank $\mathcal{H}^2$-representation which paves the way to the $O(N)$ iterative and $O(N log N)$ direct VIE solvers to be described in next section. To help readers better understand the proposed two-stage algorithm, it is necessary to first review the cluster tree structure used to model an $\mathcal{H}^2$-matrix and the $\mathcal{H}^2$-matrix partitioning.

### A. Cluster Tree and $\mathcal{H}^2$-Matrix Partitioning

A cluster tree captures the hierarchical dependence of the entire unknowns to be solved in a given problem. To build a cluster tree, we recursively split the 3-D computational domain that is composed of the SWG basis functions for $D_n$ into two sub-domains till the number of unknowns in each sub-domain becomes less than or equal to the *leafsize* ($n_{min}$). Typical values of $n_{min}$ can be as small as 2 to as large as 100. The process results in a cluster tree with each node in the tree called a cluster, as illustrated by the left (right) tree of Fig. 4(a). The root cluster is nothing but the entire unknown set, and the clusters at the bottom leaf level correspond to the subdomains whose unknown number is no greater than *leafsize*. The *leafsize* essentially controls the depth of the tree. While doing the splitting operation, special care is taken to adaptively make sure that at each nonleaf tree level, every cluster has two children of similar size. We call this splitting as a balanced splitting. Such a splitting facilitates most efficient computational cases in the arithmetic operations to be performed with the resultant matrix.

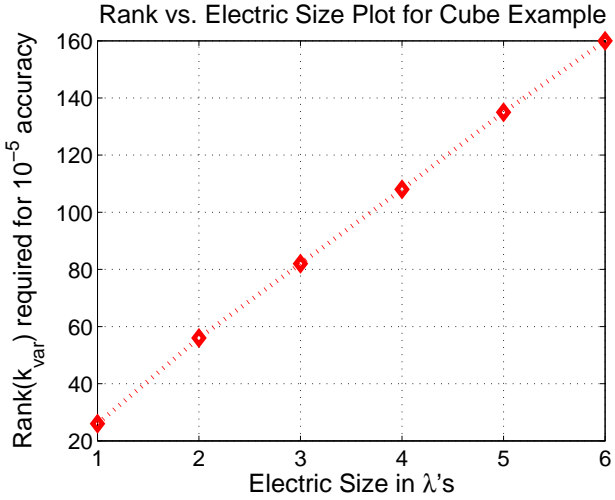

Fig. 3. Dielectric cube array (3-D structure): SVD rank growth w.r.t. the electric size in an off-diagonal block.

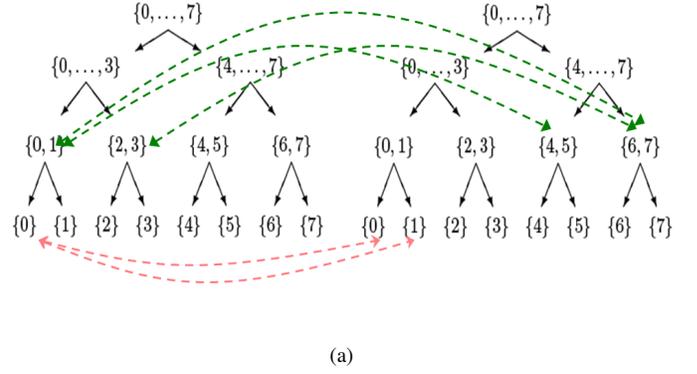

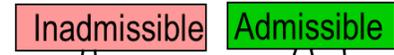

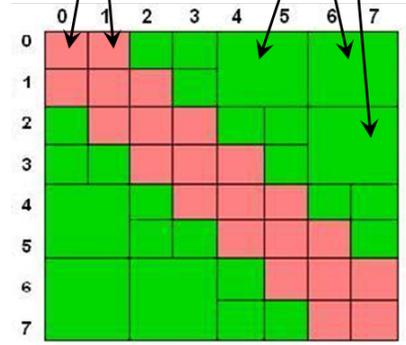

(a)

(b)

Fig. 4. Illustration of a block cluster tree and resulting $\mathcal{H}^2$-matrix partition. (a) Block cluster tree. (b) $\mathcal{H}^2$-matrix structure.

From the cluster tree, we partition the original dense VIE matrix into multilevel admissible blocks based on a strong admissibility condition [42]. To explain the process, we can place the cluster tree in parallel with itself as illustrated in Fig. 4(a). We can call the left tree row tree, and the right one column tree, as their interaction forms a matrix. Starting from the root level, we level-by-level check whether a cluster in the row tree (denoted by $t$), and a cluster in the column tree (denoted by $s$) are admissible or not. The two clusters are said to be admissible if they satisfy the following strong admissibility condition [42]

$$\max\{diam(\Omega_t), diam(\Omega_s)\} \leq \eta dist(\Omega_t, \Omega_s), \quad (4)$$

where $\eta$ is a positive parameter, and $diam(.)$ and $dist(.,.)$ respectively denote the Euclidean diameter of the support of a cluster denoted by $\Omega$, and Euclidean distance between the supports of any two clusters. As apparent from the condition that once clusters $t$ and $s$ are admissible, they ought to be physically apart. If $t$ and $s$ are admissible, they form an admissible block at that tree level, and we do not check their children clusters. Such an admissible block is marked in green in Fig. 4(b), denoted by a green link in Fig. 4(a). If they are not admissible, we proceed to check whether their children clusters satisfy the admissibility condition or not. This procedure continues until we reach the leaf level. At the end, the original dense VIE system matrix is partitioned into multilevel admissible and inadmissible blocks, as illustrated in Fig. 4(b). All the inadmissible blocks are formed at the leaf level between leaf clusters.

### B. Stage I: Cluster-Based $\mathbf{AB}^T$-Representation

With the $\mathcal{H}^2$ cluster tree and matrix partition built, we first construct a *cluster*-based non-degenerate hierarchical low-rank representation. This is different from a commonly used block-based low-rank representation like that in an $\mathcal{H}$-matrix [47]. With this new approach, the admissible blocks formed by a single cluster at its tree level are grouped together to be represented by a single low-rank block. Such a representation significantly reduces the storage and time requirements for obtaining a low-rank representation of the original dense system matrix, especially in 3-D settings. This makes the solution of millions of unknowns resulting from the electrically large VIE feasible on a single-core CPU. The accuracy of the low-rank representation of the multiple admissible blocks formed by a single cluster is also better controlled since their weights in the matrix relative to each other are considered, as compared to individually building a low-rank form for each admissible block. After the initial cluster-based low-rank representation is generated, we proceed to the second stage to obtain a minimal-rank $\mathcal{H}^2$-matrix.

For each cluster in the cluster tree, it can form multiple admissible blocks at its tree level as shown by the links in Fig. 4(a). The number of such admissible blocks is bounded by a constant $C_{sp}$ [46]. The $C_{sp}$ can be as large as hundreds in a 3-D configuration. Hence, it is not efficient to handle each admissible block one by one and generate its low-rank form individually. Instead, we propose to group these admissible blocks together and generate a *single* low-rank representation. Although physically, these blocks can be scattered in the matrix as disconnected blocks, algorithm wise, we can put them together to form a single block. After generating the low-rank representation for this single block, we can distribute it back to the original location of each admissible block if needed. As a result, our low-rank representation has a one-to-one correspondence with each cluster in an $\mathcal{H}^2$-tree, instead of being individually constructed for each admissible block. This concept is illustrated in Fig. 5.

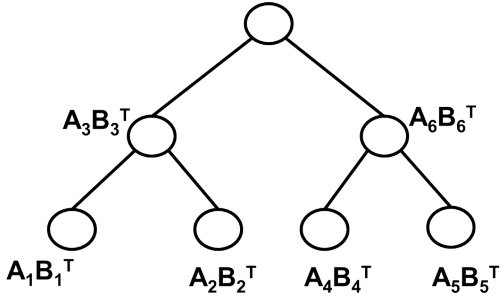

Fig. 5. Cluster-based low-rank representation.

Consider an arbitrary cluster $t$, and its associated admissible blocks $\mathbf{G}^{(t,s_1)}$, $\mathbf{G}^{(t,s_2)}$, ..., and $\mathbf{G}^{(t,s_p)}$ at its tree level. We group them into a single matrix as the following

$$\mathbf{G}^t = \left[\mathbf{G}^{(t,s_1)}, \mathbf{G}^{(t,s_2)}, ..., \mathbf{G}^{(t,s_p)}\right], \tag{5}$$

where $p$ is the number of admissible blocks formed by $t$ at its tree level. We then perform an Adaptive Cross Approximation (ACA) [ [47], pg.69] on this matrix and obtain a factorized form of

$$\mathbf{G}^t \overset{\epsilon_A}{=} \tilde{\mathbf{A}}_{\#t \times k'} \tilde{\mathbf{B}}^T_{(\#s_1+\#s_2+...+\#s_p) \times k'}, \tag{6}$$

based on prescribed accuracy $\epsilon_A$. The computational cost of this step is simply $O(k'^2(\#t+\#s))$, where $\#s = \#s_1+\#s_2+...+\#s_p$, and $\#$ denotes the cardinality of a set. This is much smaller than the cost of a brute-force SVD which scales cubically with block dimension. It is also smaller than a block-by-block ACA procedure, whose cost is $O(k'^2(\#t)C_{sp}+k^2\#s)$. In addition, the storage of $\mathbf{G}^t$ is also greatly reduced from $O(k'(\#t)C_{sp}+k\#s)$ required by a block-based ACA to $O(k'(\#t+\#s))$ units. In addition to computational efficiency, the accuracy of such a low-rank representation is also better controlled since now all admissible blocks are put together and their approximation error is controlled by $\epsilon_A$ as a whole in (6). The resulting rank is also the minimal one required to represent the entire admissible blocks formed by cluster $t$ for $\epsilon_A$ accuracy in the ACA procedure. For example, for an admissible block whose column cluster is very far from its row cluster, to achieve a $1\%$ accuracy in representing itself may require a rank of 40. However, when being put together with other admissible blocks formed by the same cluster, this block may not contribute any additional rank as it is negligible in matrix norm. Hence, it is more efficient to group the admissible blocks together to construct a single low-rank block as the relative weight of each admissible block in the entire matrix is taken into account in this representation.

After ACA, we perform another SVD, obtaining

$$\mathbf{G}^t \overset{\epsilon_{acc}}{=} \mathbf{A}_{\#t \times k} \mathbf{B}^T_{\#s \times k}, \tag{7}$$

where rank $k$ is further reduced based on the required accuracy $\epsilon_{acc}$. This step is performed because the rank $k$ determined from ACA is not the minimal rank determined by accuracy [33]. Furthermore, because the initial matrix to perform SVD is a factorized low-rank form obtained from ACA, the SVD can be performed as a reduced SVD (r-SVD) [47], whose computational cost is much reduced to $O(k^2(\#t+\#s))$ as compared to a brute-force SVD. As shown in [33], adding an additional step of r-SVD after ACA is effective in further reducing the rank without sacrificing prescribed accuracy.

After the aforementioned procedure, we obtain a cluster-based $\mathbf{AB}^T$ representation of the original VIE system matrix, as illustrated in Fig. 5. The rank of this representation is minimized based on accuracy for each cluster.

### C. Stage II: Minimal-Rank $\mathcal{H}^2$-Matrix Construction

With the $\mathbf{AB}^T$ form obtained for each cluster of an $\mathcal{H}^2$-tree, in this section, we show how to construct an $\mathcal{H}^2$-matrix representation out of it. Such a minimal-rank $\mathcal{H}^2$-matrix provides a compact nested structure for more efficient computation.

In an $\mathcal{H}^2$-matrix, the inadmissible blocks keep their original full-matrix representations, while admissible blocks are represented in the following factorized form:

$$\mathbf{G}^{t,s} = \mathbf{V}^t \mathbf{S}^{t,s} \mathbf{V}^{sT}; \tag{8}$$
$$\mathbf{V}^t \in \mathbb{C}^{\#t \times k^t}, \mathbf{S}^{t,s} \in \mathbb{C}^{k^t \times k^s}, \mathbf{V}^s \in \mathbb{C}^{\#s \times k^s}$$

with $\mathbf{V}^t(\mathbf{V}^s)$ called a cluster basis associated with $t$ and $s$ respectively. $\mathbf{S}^{t,s}$ is called a coupling matrix, and $k^{t(s)}$ is the rank of $\mathbf{V}^t(\mathbf{V}^s)$. The cluster basis $\mathbf{V}^t$ is nested in an $\mathcal{H}^2$-matrix, which satisfies the following relation

$$\mathbf{V}^t = \begin{bmatrix} \mathbf{V}^{t_1}\mathbf{T}^{t_1} \\ \mathbf{V}^{t_2}\mathbf{T}^{t_2} \end{bmatrix} = \begin{bmatrix} \mathbf{V}^{t_1} & 0 \\ 0 & \mathbf{V}^{t_2} \end{bmatrix} \begin{bmatrix} \mathbf{T}^{t_1} \\ \mathbf{T}^{t_2} \end{bmatrix} \tag{9}$$

where $t_1, t_2 \in children(t)$ are the two children clusters of $t$. $\mathbf{T}^{t_1}$ and $\mathbf{T}^{t_2}$ are called transfer matrices associated with a non-leaf cluster $t$, and they are used to build a relationship between $t$ and its two children.

It can be seen clearly that to build an $\mathcal{H}^2$-representation, we need to find cluster basis $\mathbf{V}$ for each leaf cluster, transfer matrices $\mathbf{T}^{t_1}$ and $\mathbf{T}^{t_2}$ for each non-leaf cluster, as well as coupling matrix $\mathbf{S}$ for each admissible block. In the following, we show how to find them efficiently from the $\mathbf{AB}^T$ form obtained for each cluster at the first stage. This algorithm is a bottom-up tree-traversal procedure, where we first compute the cluster basis for the leaf clusters, then obtain the two transfer matrices of each non-leaf cluster level by level. After the cluster bases are obtained, we compute the coupling matrices.

*1) Leaf Clusters:* We first generate cluster bases at the leaf level. To build a nested relationship among cluster bases, for each leaf cluster, we require its cluster basis not only represent the admissible block formed by this cluster at the leaf level, but also the admissible blocks involving this cluster at all the other levels. To do so, for each leaf cluster $t$, we build the following Gram matrix:

$$\begin{aligned}
\mathbf{G}_2^{\mathbf{t}} &= \mathbf{G}_{(t)}\mathbf{G}_{(t)}{}^H \\
&= \mathbf{A}_i(\mathbf{B}_i^T\bar{\mathbf{B}}_i)\mathbf{A}_i^H \\
&\quad + \sum_{j=ancestral-level} \mathbf{A}_{t,j}(\mathbf{B}_j^T\bar{\mathbf{B}}_j)\mathbf{A}_{t,j}^H,
\end{aligned} \tag{10}$$

where $\mathbf{G}_{(t)}$ represents the low-rank block in $\mathbf{G}$ whose rows correspond to the unknowns in cluster $t$. This block is composed of a single $\mathbf{A}_i\mathbf{B}_i{}^T$ (7) that captures all the admissible blocks formed by $t$ at $t$'s level, as well as $t$-related rows in the $\mathbf{A}_j\mathbf{B}_j{}^T$ blocks formed by $t$'s parent clusters at non-leaf levels. Here, $j$ is the index of the parent clusters of $t$ at non-leaf levels. In (10), $\mathbf{B}_j^T\bar{\mathbf{B}}_j$ is a $k \times k$ matrix, prepared a-priori for each cluster and is just referenced to each of the lower level children. $\mathbf{A}_{t,j}$, in the summation term, represents the rows corresponding to cluster $t$ of the bigger matrix $\mathbf{A}_j$ at each ancestral level. $\bar{\mathbf{B}}_j$ is the complex conjugate of $\mathbf{B}_j$. The $\mathbf{G}_2^{\mathbf{t}}$ shown in (10) is clearly of leafsize $O(n_{min})$.

Next, we perform an accuracy controlled ($\epsilon_{acc}$) Schur or SVD decomposition to get

$$\left(\mathbf{G}_2^{\mathbf{t}}\right)_{n_{min} \times n_{min}} \overset{\epsilon_{acc}}{=} \mathbf{PDP}^H, \tag{11}$$

The cluster basis for $t$ now can be obtained as

$$\mathbf{V}^t = \mathbf{P}_{n_{min},k}, \tag{12}$$

where $k$ is determined based on $\epsilon_{acc}$ when truncating (11).

*2) Non-leaf Clusters:* For clusters at a non-leaf level, if we follow the same procedure as that in the leaf level, the Gram matrix size will become increasingly large when we traverse the tree from bottom to top. Since the cluster basis generated at the leaf level $l = L$ for a leaf cluster has already taken into account upper-level blocks related to this leaf cluster, the Gram matrix formed at one level up $l = L - 1$ can be accurately projected onto the cluster bases formed at $l = L$ level. This will yield a small $k \times k$ matrix for which the cost of SVD is trivial. Similarly, the Gram matrix formed at $l = L - 2$ level can be accurately projected onto the cluster bases formed at $l = L-1$ level. Hence, the non-leaf cluster bases are generated level by level from bottom to top so that the entire computation becomes efficient. This enables performing an SVD on an $O(k)$ matrix at each level $l$.

Consider an arbitrary non-leaf cluster $t$, instead of directly building its Gram matrix as shown in (10), we project it onto its children's cluster bases to get a $k \times k$ matrix as the following:

$$\begin{aligned}
\mathbf{G}_{2,\mathbf{proj}}^{\mathbf{t}} &= \mathbf{A}_{i,small}(\mathbf{B}_i^T\bar{\mathbf{B}}_i)\mathbf{A}_{i_{small}}^H \\
&\quad + \sum_{j \in ances.} \mathbf{A}_{j_{small}}(\mathbf{B}_j^T\bar{\mathbf{B}}_j)\mathbf{A}_{j_{small}}^H
\end{aligned} \tag{13}$$

where

$$\mathbf{A}_{small} = \begin{bmatrix} \mathbf{V}^{t_1} & 0 \\ 0 & \mathbf{V}^{t_2} \end{bmatrix}^H \mathbf{A}, \tag{14}$$

and each of this multiplication costs $O(k^2(\#t))$ time requiring only $O(k^2)$ storage units. The remaining multiplications in (13) involve three $O(k)$-sized matrices requiring $O(k^3)$ operations only.

We then perform an SVD on $\mathbf{G}_{2,\mathbf{proj}}^{\mathbf{t}}$ to obtain

$$\left(\mathbf{G}_{2,\mathbf{proj}}^{\mathbf{t}}\right)_{k \times k} \overset{\epsilon_{acc}}{=} \mathbf{PDP}^H, \tag{15}$$

the cost of which is $O(k^3)$ only. The two transfer matrices of the non-leaf cluster $t$ now can be obtained as

$$\begin{bmatrix} \mathbf{T}^{t_1} \\ \mathbf{T}^{t_2} \end{bmatrix} = \begin{bmatrix} \mathbf{P}_1 \\ \mathbf{P}_2 \end{bmatrix}, \tag{16}$$

where $\mathbf{P}_1$ and $\mathbf{P}_2$ are $\mathbf{P}$'s block rows corresponding to $t$'s two children clusters $t_1$, and $t_2$ respectively.

*3) Formation of Coupling Matrices:* After the factorized $\mathbf{AB}^T$ form is generated for each cluster, for each admissible block in the $\mathcal{H}^2$-matrix, its factorized form is readily known as $\tilde{\mathbf{G}}^{(t,s)} = \mathbf{A}_{\#t \times k}\mathbf{B}_{\#s \times k}^T$. The $\mathbf{A}_{\#t \times k}$ is the same as that generated for cluster $t$, and the $\mathbf{B}_{\#s \times k}^T$ is simply the $t$-cluster-based $\mathbf{B}^T$'s columns corresponding to column cluster $s$.

To obtain the coupling matrix for each admissible block, we utilize the following relationship

$$\tilde{\mathbf{G}}^{(t,s)} = \mathbf{A}_{\#t \times k}\mathbf{B}_{\#s \times k}^T = \mathbf{V}^t\mathbf{S}\mathbf{V}^{s^T}, \tag{17}$$

and the fact that the proposed cluster bases are unitary. Hence, we have

$$(\mathbf{S})_{k \times k} = \mathbf{V}^{t^H}\mathbf{AB}^T\bar{\mathbf{V}}^s, \tag{18}$$

where $\bar{\mathbf{V}}^s$ is the complex conjugate of $\mathbf{V}^s$. The total cost of computing the coupling matrix is again $O(k^2(\#t + \#s))$. Since $\mathbf{V}$ is nested, at the nonleaf level, the cluster basis is manifested as transfer matrices.

*4) Transforming Complex-Valued Cluster Bases to Real-Valued Ones:* The nested cluster bases constructed in the proposed algorithm are complex-valued. They can be readily converted to real-valued cluster bases following the procedure described in Section III.E of [34]. As a result, the new cluster bases satisfy $\mathbf{V}^T\mathbf{V} = \mathbf{I}$ instead of the original $\mathbf{V}^H\mathbf{V} = \mathbf{I}$. This property helps make the subsequent matrix inverse more efficient, since our IE system matrix is symmetric instead of complex-conjugate symmetric. With the cluster bases updated to real-valued ones, the coupling matrices of the admissible blocks are also correspondingly updated as shown in Section III.E of [34].

## IV. Proposed $O(N)$ Iterative and $O(N\log N)$ Direct Volume Integral Equation Solvers for Electrically Large Analysis

### A. Storage and Complexity

The cluster bases are stored at the leaf level. For each non-leaf cluster, we store its two transfer matrices, each of size $k_{var} \times k_{var}$. Coupling matrices of size $k_{var} \times k_{var}$ are stored at corresponding levels to represent the off-diagonal admissible interactions between clusters. Considering the fact that the number of admissible blocks formed by a cluster at each level is bounded by sparsity constant $C_{sp}$, the total memory cost can be evaluated as:

$$
\begin{aligned}
Memory\ Cost &= \sum_{l=0}^{L} O(k_{var}^2(2^l + C_{sp}2^l)) \\
&\quad + 2C_{sp}O(n_{min}^2)N \\
&= \sum_{l=0}^{L} O(k_{var}^2 C_{sp}2^l),
\end{aligned}
\tag{19}
$$

where the $n_{min}$-related term is associated with the storage of inadmissible blocks.

With the proposed minimal-rank $\mathcal{H}^2$-representation, the rank ($k_{var}$) scales linearly with electrical size. Hence, for VIE, we have

$$
k_{var} = O(N^{\frac{1}{3}}).
\tag{20}
$$

Substituting it into (19), we obtain

$$
\begin{aligned}
Memory\ Cost &= \sum_{l=0}^{L} O\left(C_{sp}2^l \left(\frac{N}{2^l}\right)^{\frac{2}{3}}\right) \\
&= \sum_{l=0}^{L} O(C_{sp}N^{\frac{2}{3}}2^{\frac{l}{3}}) \\
&= C_{sp}O(N),
\end{aligned}
\tag{21}
$$

which is linear regardless of electrical size.

### B. Matrix-Vector Multiplication and Its Complexity

Multiplying minimal-rank $\mathcal{H}^2$-based $\mathbf{S}$ with a vector $x$ comprises of multiplying its inadmissible blocks and admissible blocks with $x$ [21], [47], [49]. In the multiplication with admissible blocks, we can fully take advantage of the $\mathcal{H}^2$-tree structure $T_{\mathcal{I}}$ and the nested cluster bases as follows.

*1) For admissible blocks:* We perform the following three steps:
1) Forward transformation: Compute $x^s := (\mathbf{V^s})^{\mathbf{T}} x|_{\hat{s}}$ for all clusters $s \in T_{\mathcal{I}}$.
2) Coupling-matrix multiplication: Compute $y^t := \sum_{s \in R^t} \mathbf{S}^{t,s} x^s$ for all clusters $t \in T_{\mathcal{I}}$ where $R^t$ contains all the clusters $s$ such that $(t, s)$ is an admissible block.
3) Backward transformation: Compute $y_i := \sum_{t, i \in \hat{t}}(\mathbf{V^t} y^t)_i$.

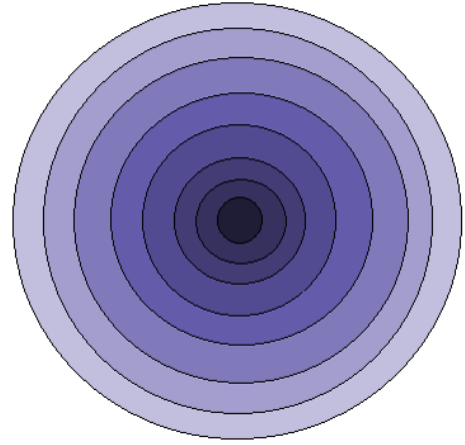

(a)

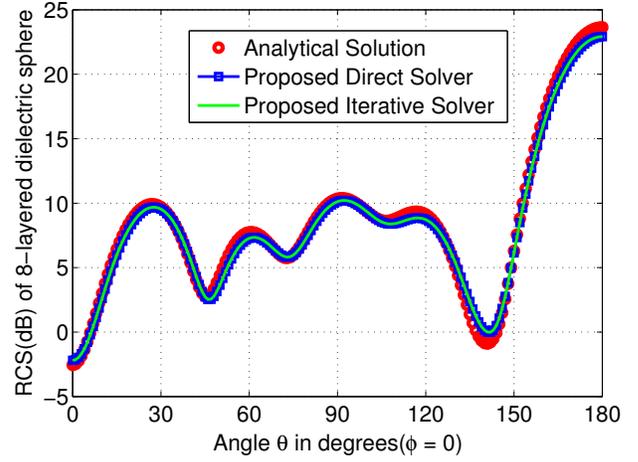

(b)

Fig. 6. Multilayer dielectric sphere of $k_0 a = 6.28$: structure and RCS comparison. (a) Multilayer dielectric sphere. (b) RCS comparison with Mie-series solution.

*2) For inadmissible blocks:* For these blocks, a full matrix-vector multiplication is performed.

With the proposed minimal-rank $\mathcal{H}^2$-representation, the total operation count for a matrix-vector multiplication is given by

$$
\begin{aligned}
MVM\ Cost &= \sum_{l=0}^{L} O(C_{sp}2^l k_{var}^2) = \sum_{l=0}^{L} O\left(C_{sp}2^l \left(\frac{N}{2^l}\right)^{\frac{2}{3}}\right) \\
&= \sum_{l=0}^{L} O(C_{sp}N^{\frac{2}{3}}2^{\frac{l}{3}}) \\
&= C_{sp}O(N).
\end{aligned}
\tag{22}
$$

We can see that with a minimal-rank $\mathcal{H}^2$-representation, even though the rank grows with electrical size, the complexity of $\mathcal{H}^2$-based matrix-vector multiplication is kept linear for general VIE-based electrodynamic analysis.

## C. Matrix Inversion and Its Complexity

The $\mathcal{H}^2$-represented system matrix $\mathbf{S}$ can be written as

$$\mathbf{S} = \begin{bmatrix} \mathbf{S_{11}} & \mathbf{S_{12}} \\ \mathbf{S_{21}} & \mathbf{S_{22}} \end{bmatrix},$$

then by Sherman-Morrison-Woodbury Formula, the inverse for $\mathbf{S}$ given by $\mathbf{S^{-1}}$ can be recursively evaluated as:

$$\mathbf{S^{-1}} = \begin{bmatrix} \mathbf{S_{11}^{-1}} \oplus \mathbf{S_{11}^{-1}} \otimes \mathbf{S_{12}} \otimes \mathbf{F^{-1}} \otimes \mathbf{S_{21}} \otimes \mathbf{S_{11}^{-1}} & -\mathbf{S_{11}^{-1}} \otimes \mathbf{S_{12}} \otimes \mathbf{F^{-1}} \\ -\mathbf{F^{-1}} \otimes \mathbf{S_{12}} \otimes \mathbf{S_{11}^{-1}} & \mathbf{F^{-1}} \end{bmatrix}$$

where $\mathbf{F} = \mathbf{S_{22}} \oplus (-\mathbf{S_{21}} \otimes \mathbf{S_{11}^{-1}} \otimes \mathbf{S_{21}})$ and $\oplus, \otimes$ are addition and multiplication defined for the $\mathcal{H}^2$-matrix elaborated. The recursive inverse formulation can be realized by the code in Table I. It can be seen that the computation of inverse

### TABLE I
### RECURSIVE INVERSE ALGORITHM

| Procedure $\mathcal{H}^2$-inverse($\mathbf{S}, \mathbf{X}$) ($\mathbf{X}$ is temporary memory) |
|---|
| If matrix $\mathbf{S}$ is a non-leaf matrix block |
| $\quad$ $\mathcal{H}^2$-inverse ($\underline{\mathbf{S_{11}}}, \mathbf{X_{11}}$) |
| $\quad$ $\mathbf{S_{21}} \otimes \underline{\mathbf{S_{11}}} \to \mathbf{X_{21}}$ |
| $\quad$ $\underline{\mathbf{S_{11}}} \otimes \mathbf{S_{21}} \to \mathbf{X_{12}}$ |
| $\quad$ $\overline{\mathbf{S_{22}}} \oplus (-\mathbf{X_{21}} \otimes \mathbf{S_{12}}) \to \underline{\mathbf{S_{22}}}$ |
| $\quad$ $\mathcal{H}^2$-inverse ($\underline{\mathbf{S_{22}}}, \mathbf{X_{22}}$) |
| $\quad$ $-\underline{\mathbf{S_{22}}} \otimes \mathbf{X_{21}} \to \mathbf{S_{21}}$ |
| $\quad$ $-\mathbf{X_{12}} \otimes \underline{\mathbf{S_{22}}} \to \mathbf{S_{12}}$ |
| $\quad$ $\underline{\mathbf{S_{11}}} \oplus (-\overline{\mathbf{S_{12}}} \otimes \mathbf{X_{21}}) \to \mathbf{S_{11}}$ |
| else |
| $\quad$ Inverse (S) (normal full matrix inverse) |

involves a full-matrix inverse at the leaf level and a number of block matrix-matrix multiplications at the non-leaf levels. The main operation in the inverse algorithm is to perform fast block matrix multiplications based on orthogonal nested cluster basis. For example, since $\mathbf{V}^T\mathbf{V} = \mathbf{I}$ is satisfied for each cluster $s$, an admissible block based matrix multiplication encountered in the inverse procedure can be done based on

$$\mathbf{V}^t\mathbf{S_1}\mathbf{V}^{sT} \times \mathbf{V}^s\mathbf{S_2}\mathbf{V}^{rT} = \mathbf{V}^t\mathbf{S_1}(\mathbf{V}^{sT} \times \mathbf{V}^s)\mathbf{S_2}\mathbf{V}^{rT}$$
$$= \mathbf{V}^t\mathbf{S_1}(\mathbf{I})\mathbf{S_2}\mathbf{V}^{rT}$$
$$= \mathbf{V}^t\mathbf{S_1}\mathbf{S_2}\mathbf{V}^{rT}$$

where only $\mathbf{S_1}\mathbf{S_2}$ needs to be computed, the cost of which is $O(k_{var}^3)$. Similarly, all the remaining 6 cases which are encountered out of the possible 27, have been presented in [33]. Each is bounded by $O(k_{var}^3)$ operations. As shown in [30], each of the $C_{sp}2^l$ number of admissible blocks at level $l$ requires $C_{sp}$ number of block multiplications each costing $O(k^3_{var})$ operations. So, the total cost to get an $\mathcal{H}^2$-matrix based inverse is:

$$Inversion\ Cost$$
$$= \sum_{l=0}^{L} (\#\ of\ blocks\ at\ level\ l)(one\ block\ cost)$$
$$= \sum_{l=0}^{L} (C_{sp}2^l)O(C_{sp}k_{var}^3)$$
$$= \sum_{l=0}^{L} C_{sp}^2 2^l O(k_{var}^3) \qquad (23)$$

With the rank's growth with electrical size taken into account, we obtain

$$Inversion\ Cost = \sum_{l=0}^{L} O\left((\frac{N}{2^l})^{\frac{3}{3}}C_{sp}^2 2^l\right)$$
$$= O\left(C_{sp}^2 N \sum_{l=0}^{L} 1\right)$$
$$= C_{sp}^2 O\left(NlogN\right) \qquad (24)$$

Thus we see that for a VIE, with the proposed minimal-rank representation, the underlying storage and matrix-vector multiplication cost becomes linear while the inversion cost becomes as fast as $O(NlogN)$.

## D. Matrix-Vector and Matrix-Matrix Multiplication for Solutions

The solution vector $D$ is then obtained by multiplying the resulting $\mathbf{S^{-1}}$ having the same $\mathcal{H}^2$-structure as $\mathbf{S}$, by $E$. For multiple right hand sides, we use $\mathcal{H}^2$-based matrix-matrix multiplication to obtain final solutions.

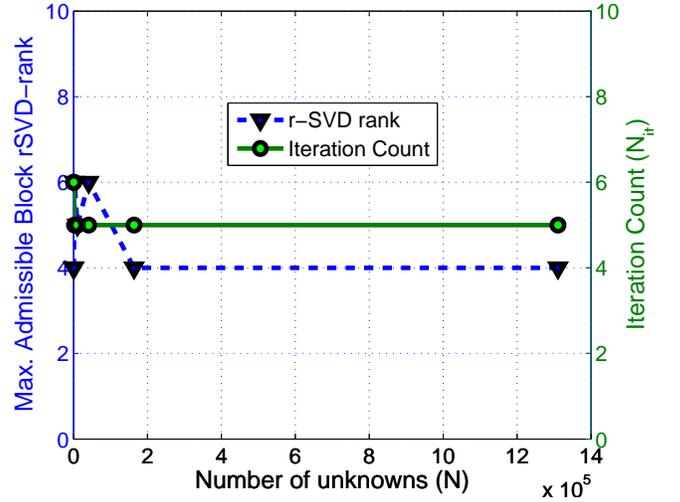

Fig. 7. Rank and iteration counts.

## V. NUMERICAL RESULTS

The numerical results from the proposed fast solvers are validated with the analytical Mie Series solution in the first multi-layered dielectric sphere examples. Next, large-scale 1-, 2- and 3-D dielectric structures, resulting in more than millions of unknowns, are simulated to demonstrate the accuracy controlled performance benefits that can be achieved with the proposed solvers. In all these numerical examples, $\eta = 1$ is used in (4) for building the $\mathcal{H}^2$ cluster tree. The computer used has a single Inte Xeon E5-2690 CPU core running at 3 GHz.

### A. Analytical Validation for an Eight Layered Dielectric Sphere

The numerical results from our proposed fast solvers are first validated with the analytical Mie Series. An eight layered dielectric sphere of 1 wavelength (free-space) radius is simulated with a permittivity profile which increases from the outermost to the center of the sphere. The outermost layer has a relative permittivity of 1.5 which increases in steps of 0.5 each layer to the innermost layer value of 5.0. The cross-sectional view of the simulated sphere is shown in Fig. 6(a) . The field of excitation is a normalized $-z$ directed plane wave polarized along the $x$-axis of standard cartesian coordinates. In Fig. 6(b), radar-cross-section (RCS in dB) of this eight layered sphere is plotted as a function of spherical coordinate polar angle ($\theta$ in degrees) for zero azimuth. It is evident that the numerical results from our proposed solvers (both direct and iterative solvers) show good agreement with the analytical Mie Series solution. The convergence criterion used in the BiCGStab iterative solver is $10^{-3}$.

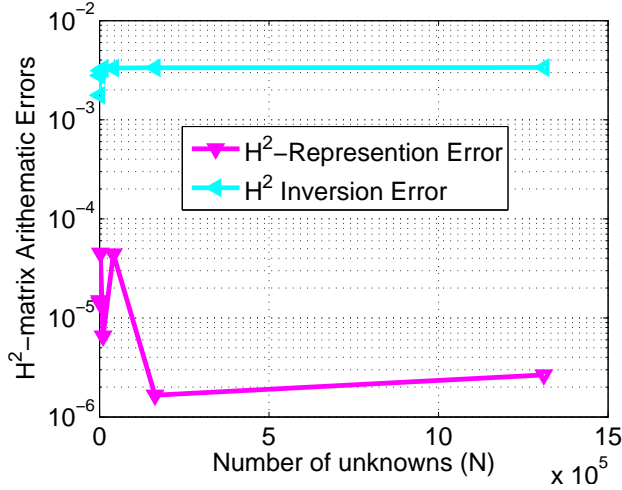

Fig. 8. $\mathcal{H}^2$-matrix representation and inversion accuracy.

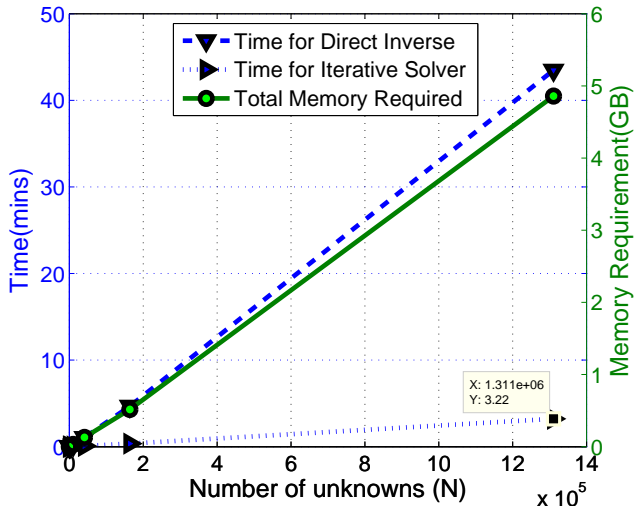

Fig. 9. Dielectric Rod: Scaling of computational resources with number of unknowns.

### B. Large-scale Dielectric Rod

The performance of the proposed solvers is first demonstrated on a dielectric rod of relative permittivity 2.54. The cross-section of the rod is fixed at $\lambda_0/10 \times \lambda_0/10$ whereas its length is increased from $1\lambda_0$ to $8,194\lambda_0$. This results in scaling the number of unknowns, $N$, from 164 to as large as 1.31 million. A normalized $-y$ directed and $z$-polarized plane wave is used as excitation. Since only one dimension of the simulated structure is changing, this example essentially represents 1-D problem characteristics [51]. In Fig. 7, it is proven numerically that indeed the accuracy determined rank for all these simulated large-scale rod examples remain constant, irrespective of the electric size. We also plot the total number of iterations required for BiCGStab iterative solver for a convergence criterion of $10^{-3}$. It is clear that only 5 iterations are required to converge regardless of the electric size of the rod. In Fig. 8, the representation and inversion accuracy is shown as a function of the scaling of number of unknowns. It is clear that the accuracy of the inverse is well-controlled below 0.5% error for all the simulated large-scale problems.

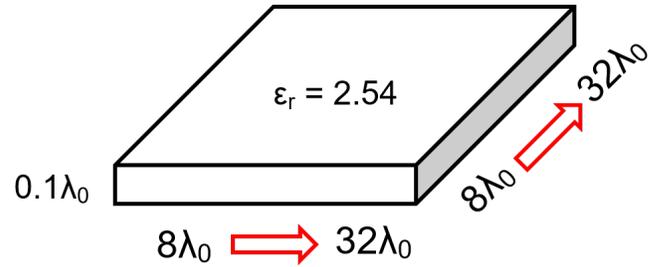

Fig. 10. Simulated dielectric slab.

Finally, the scaling of computational resources with number of unknowns is presented in Fig. 9. It is evident from the numerical results that, because of constant representation rank, storage requirements and solution times for iterative and direct solvers scale linearly with the number of unknowns. This proves the validity of equations (22) and (24). We can see that even for a dense matrix of size 1.31 million unknowns, it only takes about 40 minutes to get an accuracy controlled inverse of the VIE system matrix. For such rod-like structures, our proposed fast iterative solution can be obtained in less than 4 minutes.

### C. Large-scale Dielectric Slab

For 2-D field variation, the performance of the proposed solvers is demonstrated on a dielectric slab of relative permittivity 2.54. The thickness of the slab is fixed at $\lambda_0/10$ whereas its length and width is increased from $8\lambda_0 \times 8\lambda_0$ to $16\lambda_0 \times 16\lambda_0$ and finally to $32\lambda_0 \times 32\lambda_0$. Such a dimension scaling results in the scaling of the number of unknowns, $N$, from 89,920 to as large as 1,434,880 that is over 1.43 million. The geometry of the simulated structure is shown in Fig. 10. Normalized $-y$ directed and $z$-polarized plane wave is used as excitation. Since two dimensions of the simulated structure are changing, this example represents 2-D problem characteristics [51]. In Fig. 11, it is proven numerically that for maintaining the same level of matrix representation accuracy ($\sim 0.03\%$), rank for these simulated large-scale slab structures increases in a square-root of logarithmic trend with the electric size. The performance of the iterative solver is summarized in Fig. 12(a). A single matrix-vector multiplication time scales linearly with the number of unknowns as well as the memory consumption. In Fig. 12(b), the direct inverse time is shown as a function of $N$. It is clear that the inverse time scales linearly with the number of unknowns since the rank scales as square-root of log-linear of electrical size. In this figure, one more data point having $N = 876,000$ for a $25\lambda_0$-size slab is added to confirm the complexity. We have also assessed the inverse accuracy by evaluating $||\mathbf{I} - \mathbf{SS}^{-1}||$ for matrix size whose memory cost is feasible on our computer platforms. This error is shown to be 3.9e-2, and 7.49e-2 respetively for the $8\lambda_0 \times 8\lambda_0$ ($N$=89,920), and $16\lambda_0 \times 16\lambda_0$ ($N$=359,040) slab respectively.

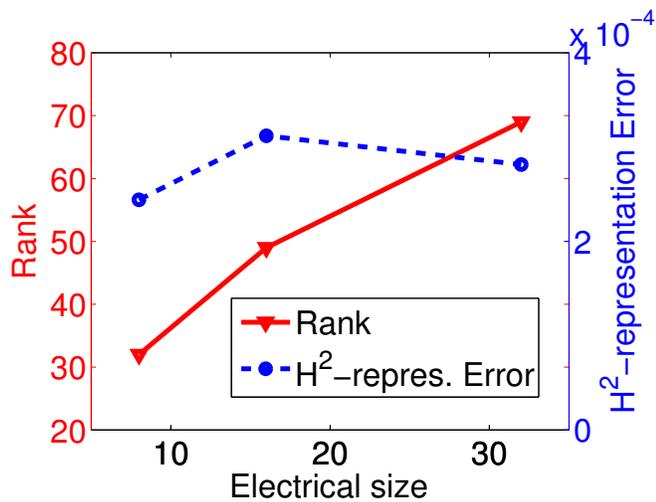

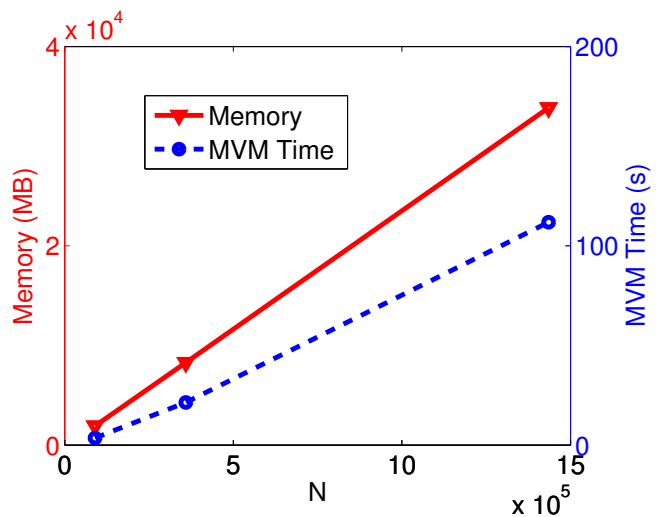

(a)

Fig. 11. Dielectric Slab: Rank and $\mathcal{H}^2$-representation error v.s. electrical size.

### D. Large-scale Array of Dielectric Cubes

A large-scale array of dielectric cubes having $\epsilon_r = 2.54$, illuminated by a plane wave ($\vec{E} = E_0 e^{jky}\hat{z}$), is simulated to demonstrate the performance of the proposed fast VIE solvers for pure 3-D field variation. The dimension of each cube is $0.3\lambda_0 \times 0.3\lambda_0 \times 0.3\lambda_0$ and the distance between neighboring cubes is fixed at $0.3\lambda_0$. The number of simulation unknowns, $N$, are scaled from $3,024$ to $1.037$ million by increasing the array size from $2 \times 2 \times 2$, $4 \times 4 \times 4$, $8 \times 8 \times 8$ to $14 \times 14 \times 14$. The geometry of the simulated structure is shown in Fig. 13. Before presenting the solvers performance, it is worth pointing out that the theoretical bounds presented in section V inherently assume that the sparsity constant, $C_{sp}$ has saturated for all tree levels. For cubic growth of unknowns for 3-D problems, as presented here, such saturation is attained in the order of millions of unknowns as shown in Table II. It is thus, important to analyze the performances for iterative and direct solvers as (Memory or MVM cost)$/C_{sp}$ and (Inverse time)$/C_{sp}^2$ respectively. In Fig. 14, it is proven numerically that for maintaining a matrix representation accuracy ($< 0.8\%$), the rank's growth rate is no greater than linear with the electric size of the cube array structure. The performance of the iterative solver follows the same trend set by rank representation and is presented in Fig. 15(a). A single matrix-vector multiplication time as well as storage scales linearly with the number of unknowns. In Fig. 15(b), the direct inverse time divided by sparsity constant square, is plotted to show that indeed the inversion time complexity scales almost as $O(N)$, which agrees with our theoretical complexity analysis. The inverse error measured by $||\mathbf{I} - \mathbf{SS}^{-1}||$ is shown to be 9.03e-3, 1.73e-2, and 3.03e-2 respectively for the array size of $2 \times 2 \times 2$, $4 \times 4 \times 4$, and $8 \times 8 \times 8$ respectively.

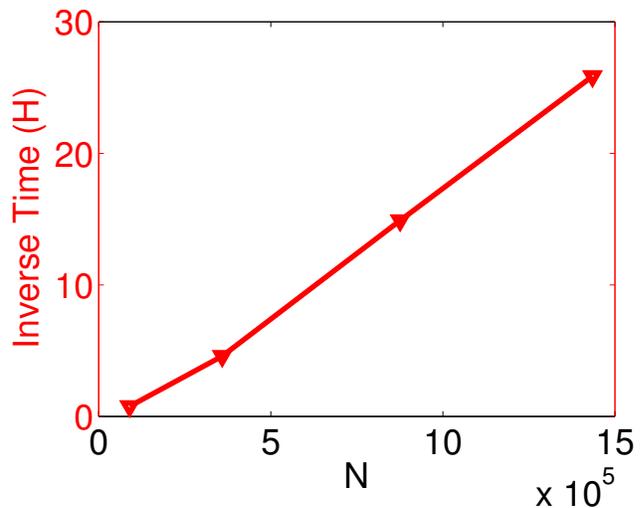

(b)

Fig. 12. Solver performance for large-scale slab structures. (a) Memory and MVM time as a function of $N$. (b) Inverse time as a function of $N$.

TABLE II
$C_{sp}$ AS A FUNCTION OF $N$ FOR THE DIELECTRIC CUBE ARRAY.

| $N$ | 3,024 | 24,192 | 193,536 | 378,000 | 1,037,232 |
|---|---|---|---|---|---|
| $C_{sp}$ | 16 | 42 | 95 | 327 | 270 |

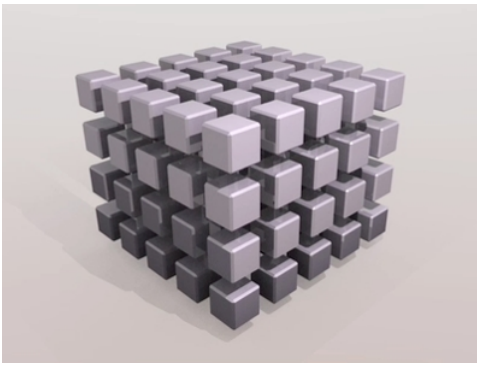

Fig. 13. Simulated arrays of large-scale dielectric cube.

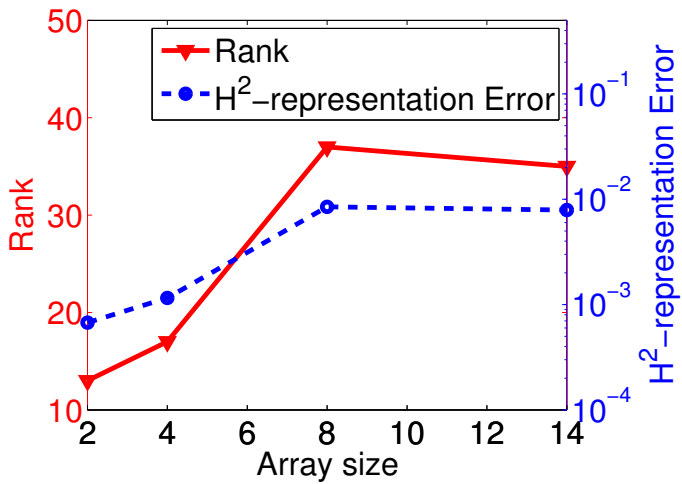

Fig. 14. Dielectric Cube Array: Rank and $\mathcal{H}^2$-representation error v.s. electrical size.

## VI. Conclusions

$O(N)$ iterative and $O(N \log N)$ direct solvers are developed for the volume integral equation based general large-scale electrodynamic analysis. The large-scale system matrix is represented by an $\mathcal{H}^2$-structure. To efficiently obtain a minimal-rank $\mathcal{H}^2$-representation, we first generate a rank-minimized $\mathbf{AB}^T$ form for each cluster in the $\mathcal{H}^2$-tree. Based on such a cluster-based low-rank form, we extract nested cluster bases and coupling matrices to obtain an efficient $\mathcal{H}^2$-matrix representation with its rank minimized by accuracy. This algorithm can also be applied to other initial representations of the IE operators such as an FMM-based representation to obtain a minimal-rank $\mathcal{H}^2$-matrix.

Analytical expressions of complexities for storage, matrix-vector-multiplication and matrix inversion are derived for general 3-D VIE-based electrodynamic analysis, confirming $O(N)$ iterative and $O(N \log N)$ direct inverse solvers with the proposed representation of VIE operators. Numerical simulations for large-scale 1-, 2- and 3-D structures, resulting in millions of unknowns, demonstrate the efficiency and complexity of the proposed VIE electrodynamic solvers. The algorithms developed in this work are kernel-independent, and hence applicable to other IE operators as well.

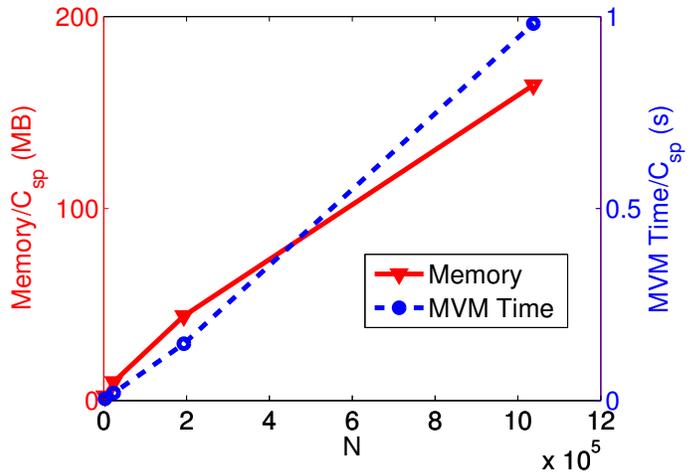

(a)

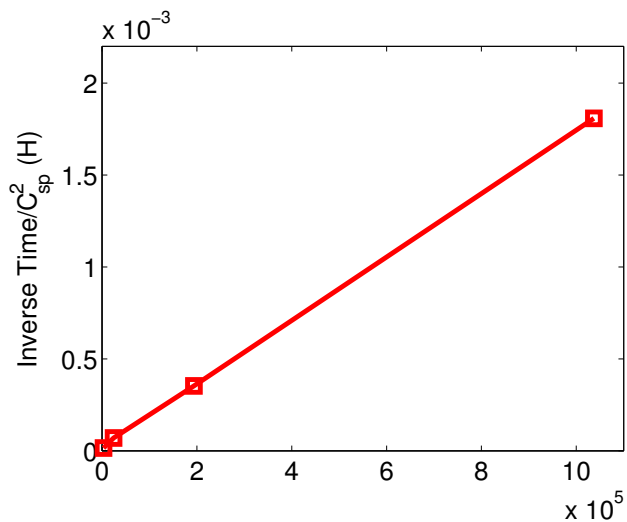

(b)

Fig. 15. Solver performance for large-scale cube array structures. (a) Memory and MVM time as a function of $N$. (b) Inverse time as a function of $N$.

## References

[1] D. H. Schaubert, D. R. Wilton, A. W. Glisson, "A tetrahedral modeling method for electromagnetic scattering by arbitrarily shaped inhomogeneous dielectric bodies," *IEEE Trans. Antennas Propag.*, vol. 32, pp. 77-85, Jan. 1984.

[2] M. M. Botha, "Solving the volume integral equations of electromagnetic scattering," *J. Comput. Phys.*, vol. 218, no. 1, pp. 141-158, 2006.

[3] M. I. Sancer, K. Sertel, J. L. Volakis, P. Van Alstine, "On volume integral equations," *IEEE Trans. Antennas Propag.*, vol. 54, pp. 1488-1495, May 2006.

[4] L. E. Sun and W. C. Chew, "A novel formulation of the volume integral equation for general large electromagnetic scattering problems," *IEEE Int. Symp. Antennas Propag.*, San Diego, CA, 2008, 4 pages.

[5] C.-C. Lu, P. Yla-Oijala, M. Taskinen, and J. Sarvas,"Comparison of two volume integral equation formulations for solving electromagnetic scattering by inhomogeneous dielectric objects," *IEEE Int. Symp. Antennas Propag.*, Charleston, SC, 2009, 4 pages.

[6] C. Pelletti, G. Bianconi, R. Mittra and A. Monorchio, "Volume integral equation analysis of thin dielectric sheet using sinusoidal macro-basis functions," *IEEE Antennas Wireless Propag. Lett.*, vol. 12, pp. 441-444, July 2009.

[7] P. Zwamborn and P. M. van den Berg, "The three-dimensional weak form

of conjugate gradient FFT method for solving scattering problems," *IEEE Trans. Microwave Theory Tech.*, vol. 40, pp. 1757-1766, Sept. 1992.

[8] H. Gan and W. C. Chew, "A discrete BCG-FFT algorithm for solving 3D inhomogeneous scatterer problems," *J. Electromag. Waves Applicat.*, vol. 9, no. 10, pp. 1339-1357, 1995.

[9] X. C. Nie, L. W. Li, N. Yuan, T. S. Yeo and Y.B. Gan, "Precorrected-FFT solution of the volume integral equation for 3-D inhomogeneous dielectric objects," *IEEE Trans. Antennas Propagat.*, vol. 53, pp. 313-320, Jan. 2005

[10] Zhenhai Zhu, Ben Song and Jacob White, "Algorithms in FastImp: A fast and wideband impedance extraction program for complicated 3-D geometries," *40th ACM/EDAC/IEEE Design Automat. Conf.*, pp. 712-717, 2003.

[11] N. A. Ozdemir, "Method of Moments solution of a non-conformal volume integral equation via the IE-FFT algorithm for electromagnetic scattering from penetrable objects," Electronic Thesis or Dissertation. Ohio State University, 2007. https://etd.ohiolink.edu/.

[12] E. Bleszynski, M. Bleszynski and T. Jarozewicz, "AIM: Adaptive integral method for solving large-scale electromagnetic scattering and radiation problems," *Radio Sci.*, vol. 31, pp. 1225-1251, Sep.-Oct. 1996.

[13] Z. Q. Zhang and Q. H. Liu, "A volume adaptive integral method (VAIM) for 3-D inhomogeneous objects," *IEEE Antennas and Wireless Propagation Letters*, vol. 1, pp. 102-105, 2002.

[14] S. Kurz, O. Rain and S. Rjasanow, "The adaptive cross-approximation technique for the 3-D boundary-element method," *IEEE Trans. Magn.*, vol. 38, pp. 421-424, 2002.

[15] N. A. Ozdemir and J.-F. Lee, "A low-rank IE-QR algorithm for matrix compression in volume integral equations," *IEEE Trans. Magn. ,* vol. 40, pp. 1017-1020, March 2004.

[16] W. C. Chew, J. M. Jin, E. Michielssen and J. M. Song, *Fast and efficient algorithms in computational electromagnetics.* Norwood, MA: Artech House, 2001.

[17] M. Kamon, M. J. Tsuk, J. K. White, "FASTHENRY: A multipole-accelerated 3-D inductance extraction program," *IEEE Trans. Microw. Theory Techn.*, vol. 42, pp. 1750-1758, Sept. 1994.

[18] K. Sertel and J. L. Volakis, "Multilevel fast multipole method solution of volume integral equations using parametric geometry modeling," *IEEE Trans. Antennas Propag.*, vol. 52, pp. 1686-1692, July 2004.

[19] C. C. Lu, "A fast algorithm based on volume integral equation for analysis of arbitrarily shaped radomes," *IEEE Trans. Antennas Propagat.*, vol. 51, no. 3, March 2003.

[20] J. Song, C. Lu, and W. Chew, "Multilevel fast multipole algorithm for electromagnetic scattering by large complex objects" *IEEE Trans. Antennas Propag.*, vol. 45, no.10, pp. 1488-1493, Oct. 1997.

[21] W. Chai and D. Jiao, "An $\mathcal{H}^2$-matrix-based fast integral-equation solver of reduced complexity and controlled accuracy for solving electrodynamic problems," *IEEE Trans. Antennas Propagat.*, vol. 57, no. 10, pp. 3147-3159, October, 2009.

[22] W. Chai and D. Jiao, "$\mathcal{H}$- and $\mathcal{H}^2$-matrix-based fast integral-equation solvers for large scale electromagnetic analysis," *IET Microwaves, Antennas & Propagation,* 4(10): pp. 1583-1596, October, 2010.

[23] E. Michielssen, A. Boag, and W. Chew, "Scattering from elongated objects: direct solution in $O(Nlog^2N)$ operations," *IEE P-Microw. Anten. P.*, pp. 277-283, 1996.

[24] P. G. Martinsson and V. Rokhlin, "A fast direct solver for scattering problems involving elongated structures," *J. Comput. Phys.*, vol. 221, pp. 288-302, 2007.

[25] R. J. Adams, Y. Xu, X. Xu, J. Choi, S. D. Gedney, and F. X. Canning "Modular fast direct electromagnetic analysis using local-global solution modes," *IEEE Trans. Antennas Propag.*, vol. 56, no. 8, Aug. 2008.

[26] J. Shaeffer, "Direct solve of electrically large integral equations for problem sizes to 1 m unknowns," *IEEE Trans. Antennas Propag.*, vol. 56, no. 8, pp. 2306-2313, Aug. 2008.

[27] E. Winebrand and A. Boag, "A multilevel fast direct solver for EM scattering from quasi-planar objects," *Proc. Int. Conf. Electromagn. Adv. Appl.*, pp. 640-643, 2009.

[28] L. Greengard, D. Gueyffier, P.-G. Martinnson,and V. Rokhlin, "Fast direct solvers for integral equations in complex three-dimensional domains," *Acta Numerica*, vol. 18, pp. 243-275, May 2009.

[29] W. Chai, D. Jiao, and C. C. Koh, "A direct integral-equation solver of linear complexity for large-scale 3D capacitance and impedance extraction," *Proc. 46th ACM/EDAC/IEEE Design Automat. Conf.*, Jul. 2009, pp. 752-757.

[30] W. Chai and D. Jiao, "Dense matrix inversion of linear complexity for integral-equation based large-scale 3-D capacitance extraction," *IEEE Trans. Microw. Theory Techn.*, vol. 59, no. 10, pp. 2404-2421, Oct. 2011.

[31] W. Chai and D. Jiao, "An LU decomposition based direct integral equation solver of linear complexity and higher-order accuracy for large-scale interconnect extraction," *IEEE Trans. Adv. Packag.*, vol. 33, no. 4, pp. 794-803, Nov. 2010.

[32] W. Chai and D. Jiao, "Direct matrix solution of linear complexity for surface integral-equation based impedance extraction of complicated 3-D structures," *Proceedings of the IEEE*, special issue on Large Scale Electromagnetic Computation for Modeling and Applications, vol. 101, pp. 372-388, Feb. 2013.

[33] W. Chai and D. Jiao, "A complexity-reduced $\mathcal{H}$-matrix based direct integral equation solver with prescribed accuracy for large-scale electro-dynamic analysis," *Proc. IEEE Int.Symp. Antennas Propag.*, Jun. 2010, DOI: 10.1109/APS.2010.5561966.

[34] W. Chai and D. Jiao, "Linear-complexity direct and iterative integral equation solvers accelerated by a new rank-minimized-representation for large-scale 3-D interconnect extraction," *IEEE Trans. Microw. Theory Techn.*, vol. 61, no. 8, pp. 2792-2805, Aug. 2013.

[35] A. Heldring, J. M. Rius, J. M. Tamayo, J. Parron, and E. Ubeda, "Multiscale compressed block decomposition for fast direct solution of method of moments linear system," *IEEE Trans. Antennas Propag.*, vol.59, no.2, pp. 526-536, Feb. 2011.

[36] A. Freni, P. D. Vita, P. Pirinoli, L. Matekovits, and G. Vecchi, "Fast-factorization acceleration of MoM compressive domain-decomposition," *IEEE Trans. Antennas Propag.*, vol.59, no.12, pp. 4588-4599, Dec. 2011.

[37] D. Gonzalez-Ovejero, F. Mesa, and C. Craeye, "Accelerated macro basis functions analysis of finite printed antenna arrays through 2D and 3D multipole expansions," *IEEE Trans. Antennas Propag.*, vol.61, no.2, pp. 707-717, Feb. 2013.

[38] Y. Brick, V. Lomakin, and A. Boag, "Fast direct solver for essentially convex scatterers using multilevel non-uniform grids," *IEEE Trans. Antennas Propag.*, vol. 62, pp. 4314-4324, 2014.

[39] S. Omar and D. Jiao, "An $\mathcal{H}^2$-Matrix based fast direct volume integral equation solver for electrodynamic analysis," *Int. Annu. Review Progress Applied Computational Electromagnetics (ACES),* Columbus, OH, 2012, 6 pages.

[40] S. Omar and D. Jiao, "A linear complexity direct volume integral equation solver for full-wave 3-D circuit extraction in inhomogeneous materials," *IEEE Trans. Microw. Theory Techn.*, vol. 63, no. 3, pp. 897-912, Mar. 2015.

[41] H. Guo, Y. Liu, J. Hu, and E. Michielssen, "A Butterfly-Based Direct Integral Equation Solver Using Hierarchical LU Factorization for Analyzing Scattering from Electrically Large Conducting Objects," *arXiv preprint arXiv:1610.00042*.

[42] S. Borm, L. Grasedyck, and W. Hackbusch, "Hierarchical Matrices," *Lecture Note 21 of the Max Planck Institute for Mathematics in the Sciences*, 2003.

[43] W. Hackbusch and B. Khoromskij, "A sparse matrix arithmetic based on $\mathcal{H}$-matrices. Part I: Introduction to $\mathcal{H}$-matrices," *Computing* , vol. 62, pp. 89-108, 1999.

[44] W. Hackbusch and B. Khoromskij, "A sparse matrix arithmetic. Part II: Application to multi-dimensional problems," *Computing*, vol. 64, pp. 21-47, 2000.

[45] S. Borm, "Introduction to hierarchical matrices with applications," *Eng. Anal. Boundary Elements (EABE)*, vol. 27, pp. 405-422, 2003.

[46] S. Borm, *Efficient Numerical Methods for Non-local Operators: $\mathcal{H}^2$-Matrix Compression Algorithms, Analysis.* Zurich, Switzerland: European Math. Soc., 2010.

[47] S. Borm, L. Grasedyck and W. Hackbusch, *"Hierarchical matrices,"* *Lecture Note 21 of the Max Plack Institute for Mathematics in the Sciences*, 2003.

[48] S. Borm, "$\mathcal{H}^2$-matrices multilevel methods for the approximation of integral operators," *Comput. Vis. Sci.*, vol. 7, pp. 173-181, 2004.

[49] S. Borm and W. Hackbusch, "$\mathcal{H}^2$-matrix approximation of integral operators by interpolation," *Appl. Numer. Math.*, vol. 43, pp. 129-143, 2002.

[50] S. Borm, "$\mathcal{H}^2$-matrix arithmetics in linear complexity," *Computing*, vol. 77, pp. 1-28, 2006.

[51] W. Chai and D. Jiao, "Theoretical study on the rank of integral operators for broadband electromagnetic modeling from static to electrodynamic frequencies," *IEEE Trans. Compon., Packag., and Manuf. Technol.*, vol. 3, no. 12, pp. 2113-2126, Dec. 2013.

[52] V. Rokhlin, "Diagonal forms of translation operators for the Helmholtz equation in three dimensions," *Appl. Comput. Harmon. Anal.*, vol. 1, pp. 8293, Dec. 1993.

[53] A. Heldring, J. Tamayo, and J. Rius, "On the degrees of freedom in the interaction between sets of elementary scatterers," *3rd European Conference on Antennas and Propagation*, 2009, pp. 2511 2514.

[54] S. Omar and D. Jiao, "An $O(N)$ iterative and $O(N log N)$ direct volume integral equation solvers for large-scale electrodynamic analysis," *the 2014 International Conference on Electromagnetics in Advanced Applications (ICEAA)*, Aug. 2014.

[55] D. Jiao and S.Omar, "Minimal-rank H2-matrix based iterative and direct volume integral equation solvers for large-scale scattering analysis," *Proc. IEEE Int. Symp. Antennas Propag.*, Jul. 2015.

[56] H. Liu and D. Jiao, "Existence of $\mathcal{H}$-matrix representations of the inverse finite-element matrix of electrodynamic problems and $\mathcal{H}$-based fast direct finite-element solvers," *IEEE Trans. Microw. Theory Techn.*, vol. 58, no. 12, pp. 3697-3709, Dec. 2010.

[57] S. Omar and D. Jiao, "An $O(N)$ Direct Volume IE Solver with a Rank-Minimized $\mathcal{H}^2$-Representation for Large-Scale 3-D Circuit Extraction in Inhomogeneous Materials," *Proc. IEEE Int. Symp. Antennas Propag.*, July 2014.